\documentclass[12pt]{article}

\usepackage[colorlinks]{hyperref}            
\usepackage{color}
\usepackage{graphicx,subfigure,amsmath,amssymb,amsfonts,bm,epsfig,epsf,url,dsfont}
\usepackage{times}
\usepackage{bbm}      
\usepackage{booktabs}
\usepackage{cases}
\usepackage{fullpage}
\usepackage[small,bf]{caption}
\usepackage{natbib}
\usepackage[top=1in,bottom=1in,left=1in,right=1in]{geometry}

\renewcommand{\phi}{\varphi}

\newcommand{\mS}{\mathbb{S}}
\renewcommand{\P}{\mathbb{P}}
\newcommand{\E}{\mathbb{E}}

\newcommand{\R}{\mathbb{R}}

\newcommand{\cK}{\mathcal{K}}

\newcommand{\cN}{\mathcal{N}}

\def\ds1{\mathds{1}}
\renewcommand{\epsilon}{\varepsilon}
\newcommand{\eps}{\epsilon}
\newcommand\Var{{\dsV\text{ar}}\,}
\newcommand\dsV{\mathbb{V}}

\newcommand{\argmin}{\mathop{\mathrm{argmin}}}

\renewcommand{\tilde}{\widetilde}

\newlength{\minipagewidth}
\newcommand{\beq}{\begin{equation}}
\newcommand{\eeq}{\end{equation}}

\newcommand{\beqa}{\begin{eqnarray}}
\newcommand{\eeqa}{\end{eqnarray}}

\newcommand{\beqan}{\begin{eqnarray*}}
\newcommand{\eeqan}{\end{eqnarray*}}

\def\ba#1\ea{\begin{align*}#1\end{align*}} 
\def\banum#1\eanum{\begin{align}#1\end{align}} 

\def\eps{\varepsilon}
\def\Proj{\mathrm{Proj}}
\def\Id{\mathrm{Id}}

\def\Cov{\mathrm{Cov}}
\def\Vol{\mathrm{Vol}}
\def\diam{\mathrm{diam}}
\def\mB{\mathrm{B}}
\def\Sph{\mS^{n-1}}
\def\RR{\mathbb{R}}
\newtheorem{theorem}{Theorem}
\newtheorem{lemma}{Lemma}

\newcommand{\BlackBox}{\rule{1.5ex}{1.5ex}}  
\newenvironment{proof}{\par\noindent{\bf Proof\ }}{\hfill\BlackBox\\[2mm]}

\begin{document}

\title{Multi-scale exploration of convex functions and bandit convex optimization}

\author{S\'ebastien Bubeck
\and 
Ronen Eldan}
\date{\today}

\maketitle

\begin{abstract}
We construct a new map from a convex function to a distribution on its domain, with the property that this distribution is a multi-scale exploration of the function. We use this map to solve a decade-old open problem in adversarial bandit convex optimization by showing that the minimax regret for this problem is $\tilde{O}(\mathrm{poly}(n) \sqrt{T})$, where $n$ is the dimension and $T$ the number of rounds. This bound is obtained by studying the dual Bayesian maximin regret via the information ratio analysis of Russo and Van Roy, and then using the multi-scale exploration to solve the Bayesian problem.
\end{abstract}

\section{Introduction}
Let $\cK \subset \R^n$ be a convex body of diameter at most $1$, and $f : \cK \rightarrow [0,+\infty)$ a non-negative convex function. Suppose we want to test whether some unknown convex function $g : \cK \rightarrow \R$ is equal to $f$, with the alternative being that $g$ takes a negative value somewhere on $\cK$. In statistical terminology the null hypothesis is 
$$\mathrm{H}_0 : g = f ,$$
and the alternative is
$$\mathrm{H}_1 : \exists \ \alpha \in \cK \; \text{such that} \; g(\alpha) < - \epsilon ,$$
where $\epsilon$ is some fixed positive number. In order to decide between the null hypothesis and the alternative one is allowed to make a single noisy measurement of $g$. That is one can choose a point $x \in \cK$ (possibly at random) and obtain $g(x) + \xi$ where $\xi$ is a zero-mean random variable independent of $x$ (say $\xi \sim \cN(0,1)$). Is there a way to choose $x$ such that the total variation distance between the observed measurement under the null and the alternative is at least (up to logarithmic terms) $\epsilon / \mathrm{poly}(n)$? Observe that without the convexity assumption on $g$ this distance is always $O(\epsilon^{n+1})$, and thus a positive answer to this question would crucially rely on convexity. We show that $\epsilon / \mathrm{poly}(n)$ is indeed attainable by constructing a distribution on $\cK$ which guarantees an exploration of the convex function $f$ at {\em every scale} simultaneously. Precisely we prove the following new result on convex functions. We denote by $c$ a universal constant whose value can change at each occurence.

\begin{theorem} \label{th:exploratory}
Let $\cK \subset \R^n$ be a convex body of diameter at most $1$. Let $f : \cK \rightarrow [0,+\infty)$ be convex and $1$-Lipschitz, and let $\epsilon > 0$. There exists a probability measure $\mu$ on $\cK$ such that the following holds true. For every $\alpha \in \cK$ and for every convex and $1$-Lipschitz function $g : \cK \rightarrow \R$ satisfying $g(\alpha) < -\eps$, one has
$$
\mu\left( \left\{x \in \cK : |f(x) - g(x)| >  \frac{c}{n^{7.5} \log(1+ n/ \eps) } \max(\epsilon, f(x)) \right\} \right) >  \frac{c}{n^{3} \log(1+n/ \eps) }.
$$
\end{theorem}

Our main application of the above result is to resolve a long-standing gap in bandit convex optimization. We refer the reader to \cite{BC12} for an introduction to bandit problems (and some of their applications). The bandit convex optimization problem can be described as the following sequential game: at each time step $t=1, \hdots, T$, a player selects an action $x_t \in \cK$, and simultaneously an adversary selects a convex (and $1$-Lipschitz) loss function $\ell_t : \cK \mapsto [0,1]$. The player's feedback is its suffered loss, $\ell_t(x_t)$. We assume that the adversary is oblivious, that is the sequence of loss functions $\ell_1,\hdots,\ell_T$ is chosen before the game starts. The player has access to external randomness, and can select her action $x_t$ based on the history $H_t = (x_s, \ell_s(x_s))_{s<t}$. The player's perfomance at the end of the game is measured through the {\em regret}:
$$R_T = \sum_{t=1}^T \ell_t(x_t) - \min_{x \in \cK} \sum_{t=1}^T \ell_t(x) ,$$
which compares her cumulative loss to the best cumulative loss she could have obtained in hindsight with a fixed action, if she had known the sequence of losses played by the adversary. A major open problem since \cite{Kle04, FKM05} is to reduce the gap between the $\sqrt{T}$-lower bound and the $T^{3/4}$-upper bound for the minimax regret of bandit convex optimization. In dimension one (i.e., $\cK = [0,1]$) this gap was closed recently in \cite{BDKP15} and our main contribution is to extend this result to higher dimensions:
\begin{theorem} \label{th:main}
There exists a player's strategy such that for any sequence of convex (and $1$-Lipschitz) losses one has
$$\E R_T \leq c \ n^{11} \log^{4}(T) \sqrt{T}  ,$$
where the expectation is with respect to the player's internal randomization.
\end{theorem}
We observe that this result also improves the state of the art regret bound for the easier situation where the losses $\ell_1, \hdots, \ell_T$ form an i.i.d. sequence. In this situation the best previous bound was obtained by \cite{agarwal2011stochastic} and is $\tilde{O}(n^{16} \sqrt{T})$.

Using Theorem \ref{th:exploratory} we prove Theorem \ref{th:main} in Section \ref{sec:main}. Theorem \ref{th:exploratory} itself is proven in Section \ref{sec:exploratory}.

\section{Proof of Theorem \ref{th:main}} \label{sec:main}
Following \cite{BDKP15} we reduce the proof of Theorem \ref{th:main} to upper bounding the {\em Bayesian maximin regret} (this reduction is simply an application of Sion's minimax theorem). In other words the sequence $(\ell_1,\hdots,\ell_T)$ is now a random variable with a distribution known to the player. Expectations are now understood with respect to both the latter distribution, and possibly the randomness in the player's strategy. We denote $\E_t$ for the expectation conditionally on the random variable $H_t$. As in \cite{BDKP15} we analyze the Bayesian maximin regret with the information theoretic approach of \cite{RR14}, which we recall in the next subsection.

\subsection{The information ratio}
Let $\bar{\cK} = \{\bar{x}_1, \hdots, \bar{x}_K\}$ be a $1/\sqrt{T}$-net of $\cK$. Note that $K \leq (4 T)^n$. We define a random variable $\bar{x}^* \in \bar{\cK}$ such that
$\sum_{t=1}^T \ell_t(\bar{x}^*) = \min_{x \in \bar{\cK}} \sum_{t=1}^T \ell_t(x)$. Using that the losses are Lipschitz one has
\begin{equation} \label{eq:approxregret}
R_T \leq \sqrt{T} + \sum_{t=1}^T (\ell_t(x_t) - \ell_t(\bar{x}^*)) .
\end{equation}
We introduce the following key quantities, for $x \in \cK$,
\begin{equation} \label{eq:rtvt}
r_t(x) = \E_t (\ell_t(x) - \ell_t(\bar{x}^*)) , \;\; \text{and} \;\; v_t(x) = \mathrm{Var}_t (\E_t(\ell_t(x) | \bar{x}^*)) .
\end{equation}
In words, conditionally on the history, $r_t(x)$ is the (approximate) expected regret of playing $x$ at time $t$, and $v_t(x)$ is a proxy for the information about $\bar{x}^*$ revealed by playing $x$ at time $t$. It will be convenient to rewrite these functions slightly more explicitly. Let $i^* \in [K]$ be the random variable such that $\bar{x}^* = \bar{x}_{i^*}$. We denote by $\alpha^*$ its distribution, which we view as a point in the $K-1$ dimensional simplex. Let $\alpha_t = \E_{t} \alpha^*$. In words $\alpha_t =(\alpha_{1,t}, \hdots, \alpha_{K,t})$ is the posterior distribution of $x^*$ at time $t$. Let $f_{i,t}, f_t : \cK \rightarrow [0,1]$, $i \in [K], t \in [T]$, be defined by, for $x \in \cK$,
$$f_t(x) = \E_t \ell_t(x), \;\; f_{i,t}(x) = \E_t ( \ell_t(x) | \bar{x}^* = \bar{x}_i ) .$$
Then one can easily see that
\begin{equation} \label{eq:rtvt2}
r_t(x) = f_t(x) - \sum_{i=1}^K \alpha_{i,t} f_{i,t}(\bar{x}_i) , \;\; \text{and} \;\; v_t(x) = \sum_{i=1}^K \alpha_{i,t} (f_t(x) - f_{i,t}(x))^2 .
\end{equation}
The main observation in \cite{RR14} is the following lemma, which gives a bound on the accumulation of information (see also [Appendix B, \cite{BDKP15}] for a short proof).
\begin{lemma} \label{lem:RR}
One always has $\E \sum_{t=1}^T v_t(x_t) \leq \frac{1}{2} \log(K)$.
\end{lemma}
An important consequence of Lemma \ref{lem:RR} is the following result which follows from an application of Cauchy-Schwarz (and \eqref{eq:approxregret}):
\begin{equation} \label{eq:RR}
\E \sum_{t=1}^T r_t(x_t) \leq \sqrt{T} + C \sum_{t=1}^T \sqrt{\E v_t(x_t)} \; \Rightarrow \; \E R_T \leq 2 \sqrt{T} + C \sqrt{\frac{T}{2} \log(K)} .
\end{equation}
In particular a strategy which obtains at each time step an information proportional to its instantaneous regret has a controlled cumulative regret:
\begin{equation} \label{eq:RR2}
\E_{t} r_t(x_t) \leq \frac{1}{\sqrt{T}} + C \sqrt{\E_{t} v_t(x_t)}, \ \forall t \in [T] \; \Rightarrow \; \E R_T \leq 2 \sqrt{T} + C \sqrt{\frac{T}{2} \log(K)} .
\end{equation}
\cite{RR14} refers to the quantity $\E_{t} r_t(x_t) / \sqrt{\E_{t} v_t(x_t)}$ as the {\em information ratio}. They show that Thompson Sampling (which plays $x_t$ at random, drawn from the distribution $\alpha_t$) satisfies $\E_{t} r_t(x_t) / \sqrt{\E_{t} v_t(x_t)} \leq K$ (without any assumptions on the loss functions $\ell_t : \cK \rightarrow [0,1]$). In \cite{BDKP15} it is shown that in dimension one (i.e., $n=1$), the latter bound can be improved using the convexity of the losses by replacing $K$ with a polylogarithmic term in $K$ (Thompson Sampling is also slightly modified). In the present paper we propose a completely different strategy, which is loosely related to the Information Directed Sampling of \cite{RR14b}. We describe and analyze our new strategy in the next subsection. 

\subsection{A two-point strategy}
We describe here a new strategy to select $x_t$, conditionally on $H_t$, and show that it satisfies a bound of the form given in \eqref{eq:RR2}. To lighten notation we drop all time subscripts, e.g. one has $r(x) = f(x) - \sum_{i=1}^K \alpha_i f_i(\bar{x}_i)$, and $v(x) = \sum_{i=1}^K \alpha_i \left(f_i(x) - f(x) \right)^2$. Our objective is to describe a random variable $X \in \cK$ which satisfies
\begin{equation} \label{eq:goal}
\E r(X) \leq \frac{1}{\sqrt{T}} + C \sqrt{ \E v(X) } ,
\end{equation}
where $C$ is polylogarithmic in $K$ (recall that $K \leq (4 T)^n$). 

Let $x^* \in \argmin_{x \in \cK} f(x)$. We translate the functions so that $f(x^*)=0$ and denote $L=\sum_{i=1}^K \alpha_i f_i(\bar{x}_i)$. If $L \geq - 1/\sqrt{T}$ then $X := x^*$ satisfies \eqref{eq:goal}, and thus in the following we assume that $L \leq - 1/\sqrt{T}$. 
\newline

\noindent
\textbf{Step 1:} We claim that there exists $\epsilon \in [|L|/2, 1]$ such that
\begin{equation} \label{eq:claim1}
\alpha\left(\left\{i \in [K] : f_i(\bar{x}_i) \leq - \epsilon \right\}\right)\geq \frac{|L|}{2 \log(2/|L|) \epsilon} .
\end{equation}
Indeed assume that \eqref{eq:claim1} is false for all $\epsilon \in [|L|/2, 1]$, and let $Y$ be a random variable such that $\P(Y=-f_i(\bar{x}_i)) = \alpha_i$, then
$$|L| = \E Y \leq |L| / 2 + \int_{|L|/2}^1 \P(Y\geq x) dx < |L| / 2 + \int_{|L|/2}^1 \frac{|L|}{2 \log(2/|L|) x} dx = |L|,$$
thus leading to a contradiction. We denote $I = \left\{i \in [K] : f_i(\bar{x}_i) \leq - \epsilon \right\}$ with $\epsilon$ satisfying \eqref{eq:claim1}.
\newline 

\noindent
\textbf{Step 2:} We show here the existence of a point $\bar{x} \in \cK$ and a set $J \subset I$ such that $\alpha(J) \geq \frac{c}{n^{3} \log(1+n/ \eps) } \alpha(I)$ and for any $i \in J$, 
\begin{equation} \label{eq:goodpoints}
|f(\bar{x}) - f_i(\bar{x})| \geq \frac{c}{n^{7.5} \log(1+n/ \eps) } \max(\epsilon, f(\bar{x})).
\end{equation} 
We say that a point is {\em good} for $f_i$ if it satisfies \eqref{eq:goodpoints}, and thus we want to prove the existence of a point $\bar{x}$ which is good for a large fraction (with respect to the posterior) of the $f_i$'s. Denote
$$A_i = \left\{x \in \cK : |f({x}) - f_i({x})| \geq \frac{c}{n^{7.5} \log(1+n/ \eps)} \max(\epsilon, f({x})) \right\},$$
and let $\mu$ be the distribution given by Theorem \ref{th:exploratory}. Then one obtains:
$$\sup_{x \in \cK} \sum_{i \in I} \alpha_i \ds1\{x \in A_i\} \geq \int_{x \in \cK} \sum_{i \in I} \alpha_i \ds1\{x \in A_i\} d\mu(x) = \sum_{i \in I} \alpha_i \mu(A_i) \geq \frac{c}{n^{3} \log(1+n/ \eps) } \alpha(I) ,$$
which clearly implies the existence of $J$ and $\bar{x}$.
\newline 

\noindent
\textbf{Step 3:} 
Let $X$ be such that $\P(X=\bar{x}) = \alpha(J)$ and $\P(X=x^*) = 1-\alpha(J)$. Then
$$\E r(X) = |L| + \alpha(J) f(\bar{x}) ,$$
and using the definition of $\bar{x}$ one easily see that:
$$\sqrt{\E v(X)} \geq \sqrt{\alpha(J) v(\bar{x})} \geq \sqrt{\alpha(J) \sum_{i \in J} \alpha_i (f_i(\bar{x})-f(\bar{x}))^2} \geq \frac{c}{n^{7.5} \log(1+n/ \eps) } \alpha(J) \max(\epsilon, f(\bar{x})) .$$
Finally, since $\alpha(J) \geq \frac{c |L|}{\eps n^3 \log^2 (1+n/\epsilon)}$, the two above displays clearly implies \eqref{eq:goal}.

\section{An exploratory distribution for convex functions} \label{sec:exploratory}
In this section we construct an exploratory distribution $\mu$ of a convex function $f$ which satisfies the conditions of Theorem \ref{th:exploratory}, thus concluding the proof of Theorem \ref{th:main}. 

\subsection{The one-dimensional case}

Since our proof of Theorem \ref{th:exploratory} will proceed by induction, our first goal is to establish the result in dimension 1. This task will be much simpler than the proof for a general dimension, but already contains some of the central ideas used in the general case. In particular, a (much simpler) multi-scale argument is used.

The main ingredient is the following lemma which is easy to verify by picture (we provide a formal proof for sake of completness).
\begin{lemma} \label{lem:onedim}
	Let $f,g:\RR \to \RR$ be two convex functions. Suppose that $f(x) \geq 0$. Let $x_0,\alpha \in \RR$ be two points satisfying $\alpha-1 < x_0 < \alpha$, and suppose that $g(\alpha) < -\eps$ for some $\eps > 0$ and that 
	\begin{equation} \label{eq:assumpDerive}
	f'(x) \geq 0, ~ \forall x>x_0.
	\end{equation}
	Let $\mu$ be a probability measure supported on $[x_0, \alpha]$ whose density with respect to the Lebesgue measure is bounded from above by some $\beta>1$. Then we have
	$$
	\mu \left ( \left \{x : |f(x) - g(x)| > \tfrac 1 4 \beta^{-1} \max(\eps, f(x)) \right \} \right  ) \geq \frac{1}{2}.
	$$
\end{lemma}

\begin{proof}
	We first argue that, without loss of generality, one may assume that $f$ attains its minimum at $x_0$. Indeed, we may clearly change $f$ as we please on the interval $(-\infty, x_0)$ without affecting the assumptions or the result of the Lemma. Using the condition \eqref{eq:assumpDerive} we may therefore make this assumption legitimate.
	
	Assume, for now, that there exists $x_1 \in [x_0,\alpha]$ for which $f(x_1) = g(x_1)$. By convexity, and since $f(x_0) \geq 0$ and $g(\alpha) < 0$, if such point exists then it is unique. Let $h(x)$ be the linear function passing through $(\alpha, g(\alpha))$ and $(x_1, f(x_1))$. By convexity of $g$, we have that $|g(x) - f(x)| \geq |h(x) - f(x)|$ for all $x \in [x_0, \alpha]$. Now, since $h(\alpha) < -\eps$ and since $\alpha<x_1+1$, we have $h'(x_0) < - (\eps + f(x_0))$. Moreover, since we know that $f(x)$ is non-decreasing in $[x_0, \alpha]$, we conclude that
	\begin{align*}
	|g(x) - f(x)| ~& \geq |h(x) - f(x)| \\
	& = |h(x) - f(x_1)| + |f(x) - f(x_1)| \\
	& = (\eps + f(x_1)) |x-x_1| + |f(x) - f(x_1)| \\
	& \geq \max(\eps, f(x)) |x-x_1|, ~~ \forall x \in [x_0, \alpha].
	\end{align*}
	It follows that
	$$
	\left \{x; ~|f(x) - g(x)| <  \tfrac 1 4 \beta^{-1} \max (\eps, f(x)) \right \} \subset I := \bigl [x_1- \tfrac 1 4 \beta^{-1}, x_1 + \tfrac 1 4 \beta^{-1} \bigr ]
	$$
	but since the density of $\mu$ is bounded by $\beta$, we have $\mu(I) \leq \tfrac 1 2$ and we're done.
	
	It remains to consider the case that $g(x) < f(x)$ for all $x \in [x_0, \alpha]$. In this case, we may define 
	$$
	\tilde g(x) = g(x) + \frac{f(x_0) - g(x_0)}{\alpha-x_0} (\alpha - x).
	$$
	Note that $\tilde g(x) \geq g(x)$ for all $x \in [x_0, \alpha]$, which implies that $|g(x)-f(x)| \geq |\tilde g(x) - f(x)|$ for all $x \in [x_0, \alpha]$. Since $\tilde g(x_0) = f(x_0)$, we may continue the proof as above, replacing the function $g$ by $\tilde g$.
\end{proof}

We are now ready to prove the one dimensional case. The proof essentially invokes the above lemma on every scale between $\eps$ and $1$.

\begin{proof}[Proof of Theorem \ref{th:exploratory}, the case $n=1$]
Let $x_0 \in \cK$ be the point where the function $f$ attains its minimum and set $d = \diam(\cK)$. Define $N = \lceil \log_2 \tfrac 1 \eps \rceil + 4$. For all $0 \leq k \leq N$, consider the interval
$$
I_k = [x_0 - d 2^{-k}, x_0 + d 2^{-k}]  \cap \cK
$$
and define the measure $\mu_k$ to be the uniform measure over the interval $I_k$. Finally, we set
$$
\mu = \frac{1}{N+2} \sum_{k=0}^N \mu_k + \frac{1}{N+2}\delta_{x_0}.
$$

Now, let $\alpha \in \cK$ and let $g(x)$ be a convex function satisfying $g(\alpha) \leq -\eps$. We would like to argue that $\mu(A) \geq \frac{1}{8 \log (1 + 1/\eps)}$ for $A = \left \{ x \in \cK: ~ |f(x) - g(x)| \geq \frac{\eps}{8} \right \}$.

Set $k = \lceil \log_{1/2} (|\alpha - x_0| / d) \rceil $. Define $Q(x) = x_0 + d 2^{-k} (x-x_0)$ and set $\tilde f(x) = f(Q(x))$, $\tilde g(x) = g(Q(x))$, $\tilde \alpha = Q^{-1} (\alpha)$ and consider the interval
$$
I = Q^{-1}(I_k) \cap \{x: (x-x_0)(\alpha - x_0) \geq 0\}
$$
It is easy to check that, by definition $I$ is an interval of length $1$, contained in the interval $[x_0,\tilde \alpha]$. Defining $\tilde \mu = \mu_I$, we have that the density of $\tilde \mu$ with respect to the Lebesgue measure is equal to $1$. An application of Lemma \ref{lem:onedim} for the functions $\tilde f, \tilde g$, the points $x_0, \tilde \alpha$ and the measure $\tilde \mu$ teaches us that
\begin{align*}
\mu_k  ( A ) = \mu_{Q^{-1}(I_k)} \left ( \left \{ x: ~\left  | \tilde f(x) - \tilde g(x) \right | \geq \frac{\eps}{8}  \right \}\right ) 
\geq \frac{1}{2 }\tilde \mu \left ( \left \{ x: ~\left  | \tilde f(x) - \tilde g(x) \right | \geq \frac{\eps}{8}  \right \}\right ) \geq \frac{1}{4}.
\end{align*}
By definition of the measure $\mu$, we have that whenever $k \leq N$, one has
$$
\mu \left ( A \right ) \geq \frac{1}{N+2} \geq \frac{1}{8 \log (1 + 1/\eps)}.
$$
Finally, if $k > N$, it means that $|\alpha - x_0| < 2^{-N} < \frac{\eps}{4}$. Since the function $g$ is $1$-Lipschitz, this implies that $g(x_0) \leq - \eps/2$ which in turn gives $|f(x_0) - g(x_0)| \geq \frac{\eps}{8}$. Consequently, $x_0 \in A$ and thus $\mu(A) \geq \mu(\{x_0\}) = \frac{1}{N+2} \geq \frac{1}{8\log (1 + 1/\eps)}$. The proof is complete.
\end{proof}

\subsection{The high-dimensional case}
We now consider the case where $n \geq 2$. For a set $\Omega \subset \R^n$ and a direction $\theta \in \R^n$ we denote $S_{\Omega, \theta} = \{x \in \Omega : |\langle x, \theta \rangle| \leq 1/4\}$, and $\mu_{\Omega}$ for the uniform measure on $\Omega$. For a distribution $\mu$ we write $\Cov(\mu) = \E_{X \sim \mu} X X^{\top}$.
\medskip

As we explain in Section \ref{sec:exploratory1} our construction iteratively applies the following lemma:

\begin{lemma} \label{lem:exploratory1}
Let $\eps > 0$, $L \in [1, 2n]$. Let $\Omega \subset \R^n$ be a convex set with $0 \in \Omega$ and $\Cov(\mu_\Omega) = \mathrm{Id}$. Let $f : \Omega \rightarrow [0,\infty]$ be a convex and $L$-Lipschitz function with $f(0) = 0$. Then there exists a measure $\mu$ on $\Omega$ and a direction $\theta \in \mS^{n-1}$ such that for all $\alpha \in \Omega \setminus S_{\Omega, \theta}$ and for every convex function $g : \Omega \rightarrow \R$ satisfying $g(\alpha) < -\eps$, one has
\begin{equation} \label{eq:mu0}
\mu\left( \left\{x \in \Omega : |f(x) - g(x)| > \frac{1}{2^{50} n^{7.5} \log(1+ n / \eps) } \max(\epsilon, f(x)) \right\} \right) > \frac{1}{16 n}.
\end{equation}
\end{lemma}

The above lemma is proven in Section \ref{sec:exploratory2}. A central ingredient in its proof is, in turn, the following Lemma, which itself is proven in Section \ref{sec:exploratory3}.

\begin{lemma} \label{lem:exploratory2}
Let $\eps > 0$, $\Omega \subset \R^n$ a convex set with $\diam(\Omega) \leq M$, and $f : \Omega \rightarrow \R_+$ a convex function. Assume that there exist $\delta \in (0, \tfrac 1 {32 n^2})$, $z \in \Omega \cap \mB(0,\tfrac 1 {16})$, $\theta \in \mS^{n-1}$ and $t >0$ such that
\begin{equation} \label{eq:jollygood}
\mu_{B(z, \delta)} \left( (\nabla f)^{-1} \left( \mB\left(t \theta, \frac{t}{16n^2}\right)  \right) \right) \geq 1/2 .
\end{equation}
Then for all $\alpha \in \Omega$ satisfying $\langle \alpha , \theta \rangle \geq \tfrac{1}{8}$ and $|\alpha| \leq 2n$ and for all convex function $g : \Omega \rightarrow \R$ satisfying $g(\alpha) < -\eps$, one has
$$\mu_{\mB(z, \delta)} \left( \left\{x \in \Omega : |f(x) - g(x)| > \frac{\delta}{2^{13} M \sqrt{n}} \max(\epsilon, f(x)) \right\} \right) >  \frac{1}{8} .$$
\end{lemma}


\subsection{From Lemma \ref{lem:exploratory1} to Theorem \ref{th:exploratory}: a multi-scale exploration} \label{sec:exploratory1}

An intermediate lemma in this argument will be the following:
\begin{lemma} \label{lem:exploratory3}
There exists a universal constant $c>0$ such that the following holds true. Let $\eps > 0$, $\Omega \subset \R^n$ a convex set with $0 \in \Omega$ and $\Cov(\mu_\Omega) = \mathrm{Id}$. Let $f : \Omega \rightarrow [0,\infty)$ be a convex and $1$-Lipschitz function. Then there exists a measure $\mu$ on $\Omega$, a point $y \in \Omega$ and a direction $\theta \in \mS^{n-1}$ such that for all $\alpha \in \Omega$ satisfying
$$
\left | \langle \alpha - y, \theta \rangle \right | \geq \frac{c \eps}{16n^{10}}
$$
and for every convex function $g : \Omega \rightarrow \R$ satisfying $g(\alpha) < -\eps$, one has
\begin{equation} \label{eq:eqlem3}
\mu\left( \left\{x \in \Omega : |f(x) - g(x)| > \frac{c}{n^{7.5} \log(1+n/ \eps) } \max(\epsilon, f(x)) \right\} \right) > \frac{ c }{ n^2 \log (1+n/\eps)}.
\end{equation}
\end{lemma}

\subsubsection{From Lemma \ref{lem:exploratory3} to Theorem \ref{th:exploratory}}

Given Lemma \ref{lem:exploratory3}, the proof of Theorem \ref{th:exploratory} is carried out by induction on the dimension. The case $n=1$ has already been resolved above. Now, suppose that the theorem is true up to dimension $n-1$, where the constant $c>0$ is the constant from Lemma \ref{lem:exploratory3}. Let $\cK \in \RR^n$ and $f$ satisfy the assumptions of the theorem. Denote $Q = \Cov(\mu_{\cK})^{-1/2}$ and define
$$
\Omega = Q(\cK), ~~ \tilde f(x) = f(Q^{-1}(x))
$$
so that $\tilde f:\Omega \to \RR$. Since $\diam(\cK) \leq 1$, we know that for all $u \in \Sph$, $\Var[\Proj_u \mu_K] \leq 1$ which implies that $\| Q^{-1} \| \leq 1$. Consequently, the function $\tilde f$ is $1$-Lipschitz. We now invoke Lemma \ref{lem:exploratory3} on $\Omega$ and $\tilde f$ which outputs a measure $\mu_1$, a point $y \in \Omega$ and a direction $\theta$. By translating $f$ and $\cK$, we can assume without loss of generality that $y=0$. Fix some linear isometry $T:\RR^{n-1} \to \theta^\perp$. Define
$$
\Omega' = T^{-1} \Proj_{\theta^\perp} \left (\Omega \cap \left  \{x: |\langle x, \theta \rangle| \leq \delta \right \} \right )
$$
where $\delta = \frac{c \eps}{16n^{10}}$ and $c$ is the universal constant from Lemma \ref{lem:exploratory3}. Since $\tilde f$ is convex, there exists $I \subset \RR \times \RR^n$ so that
\begin{equation} \label{expansion}
\tilde f(x) = \sup_{(a,y) \in I} \left (a + \langle x, y \rangle \right ), ~~ \forall x \in \Omega.
\end{equation}
We may extrapolate $\tilde f(x)$ to the domain $\RR^n$ by using the above display as a definition. We now define a function $h: \Omega' \to \RR$ by
\begin{equation} \label{eq:defh}
h(x) :=\sup_{w \in [-\delta,\delta] }  \tilde f ( T(x)+w \theta ).
\end{equation}
It is clear that $\diam(\Omega') \leq 1$. Moreover, $h$ is $1$-Lipschitz since it can be written as the supremum of $1$-Lipschitz functions. 
We can therefore use the induction hypothesis with $\Omega', h(x)$ to obtain a measure $\mu_2$ on $\Omega'$. Next, for $y \in \RR^{n-1}$, define
$$
N(y) := \left \{x \in \Omega: T^{-1}(\Proj_{\theta^\perp} x ) = y  \right \}
$$
and set
$$
\mu(W) = \frac{1}{n} \mu_1(Q(W)) + \frac{n-1}{n} \int_{\Omega'} \frac{\Vol_1(Q(W) \cap N(u))}{\Vol_1(N(u)) } d \mu_2(u)
$$
for all measurable $W \subset \RR^n$. 

Fix $\alpha \in \cK$, let $g:\cK \to \RR$ be a convex and $1$-Lipschitz function satisfying $g(\alpha) \leq -\eps$. Recall that $c$ denotes the universal constant from Lemma \ref{lem:exploratory3}. Define
$$
A = \left\{x \in \cK : |f(x) - g(x)| >  \frac{c}{n^{7.5} \log(1+n/ \eps) } \max(\epsilon, f(x)) \right\}.
$$
The proof will be concluded by showing that $\mu(A) \geq \frac{c}{n^{3} \log(1+n/ \eps) }$. \\

Define $\tilde g(x) = g(Q^{-1}(x))$ and remark that $\tilde g$ is $1$-Lipschitz. First consider the case that $\left |\left \langle Q \alpha, \theta \right \rangle \right | \geq \delta$, then by construction, we have
\begin{align*}
\mu(A) ~& \geq \frac{1}{n} \mu_1(Q(A)) \\ 
&= \frac{1}{n} \mu_1 \left ( \left \{x \in \Omega; |\tilde f(x) - \tilde g(x)| >  \frac{c}{n^{7.5} \log(1+n/ \eps) } \max(\epsilon, f(x))  \right \} \right ) \\
& \stackrel{\eqref{eq:eqlem3}}{\geq} \frac{c}{n^{3} \log(1+n/ \eps) },
\end{align*}
and we're done. 

Otherwise, we need to deal with the case that $\left |\left \langle Q \alpha, \theta \right \rangle \right | < \delta$. Define $q(x)$ to be the function obtained by replacing $\tilde f(x)$ with $ \tilde g(x)$ in equation \eqref{eq:defh} and consider the set
$$
A' = \left \{ x \in \Omega'; ~ |h(x) - q(x)| > \frac{c}{(n-1)^{7.5} \log(1+n/ \eps) } \max(\epsilon, h(x)) \right\}.
$$
By construction of the measure $\mu_2$ we have $\mu_2(A') \geq \frac{c}{(n-1)^{3} \log(1+n/ \eps)}$. We claim that $N(A') \subset Q(A)$, which implies that 
$$
\mu(A) \geq \frac{n-1}{n} \mu_2(A') \geq \frac{c}{n^{3} \log(1+n/ \eps)}
$$ 
which will complete the proof. Indeed, let $y \in N(A')$. Define $z=T^{-1}(\Proj_{\theta^\perp} y)$, so that $z \in A'$. Let $w_1, w_2 \in N(z)$ be points such that
$$
h(z) = \tilde f(w_1), ~~ q(z) = \tilde g(w_2).
$$
Such points exist since, by continuity, the maximum in equation \eqref{eq:defh} is attained. Now, since $z \in A'$, we have by definition that
$$
|\tilde f(w_1) - \tilde g(w_2)| >  \frac{c}{(n-1)^{7.5} \log(1+n/ \eps) } \max(\epsilon, \tilde f(w_1)).
$$
Finally, since the functions $\tilde f,\tilde g$ are $1$-Lipschitz, we have that
\begin{align*}
|\tilde f(y) - \tilde g(y)| ~& \geq |\tilde f(w_1) - \tilde g(w_2)| - |\tilde f(y) - \tilde f(w_1)| - |\tilde g(y) - \tilde g(w_2)| \\
& \geq \frac{c}{(n-1)^{7.5} \log(1+n/ \eps) } \max(\epsilon, \tilde f(w_1)) - |y-w_1| - |y-w_2| \\
& \geq  \frac{c}{(n-1)^{7.5} \log(1+n/ \eps) } \max(\epsilon, \tilde f(y)) - 4 \delta \\
& = \frac{c}{(n-1)^{7.5} \log(1+n/ \eps) } \max(\epsilon, \tilde f(y)) - \frac{c \eps}{4 n^{10}} \\
& \geq \frac{c}{n^{7.5} \log(1+n/ \eps) } \max(\epsilon, \tilde f(y))
\end{align*}
which implies, by definition, that $y \in Q(A)$. The proof is complete.

\subsubsection{From Lemma \ref{lem:exploratory1} to Lemma \ref{lem:exploratory3}}

We construct below a decreasing sequence of domains $\Omega_0 \supset \Omega_1 \supset ... \supset \Omega_N$. Let $x_0 \in \Omega$ be a point where $f(x)$ attains its minimum on $\Omega$. Set $\Omega_0 = \Omega - x_0$. Given $i \geq 0$, we define the domain $\Omega_{i+1}$, given the domain $\Omega_i$, by induction as follows.  Define $Q_i = \Cov(\mu_{\Omega_i})^{-1/2}$ and $f_i(x) = f(Q_i^{-1} (x + x_0)) - f(x_0)$. We have
$$
|\nabla f_i(x)| = \left | Q_i^{-1} \nabla f(Q_i^{-1} (x)) \right |  \leq \| Q_i^{-1} \|.
$$
Now, by Lemma \ref{lem:tech1} we know that 
$$
\diam(\Omega_i) \leq \diam(\Omega) \leq n+1
$$
which implies that $\left \| Q_i^{-1} \right \| \leq n+1$. We conclude that $f_i$ is $(n+1)$-Lipschitz. We may therefore invoke Lemma \ref{lem:exploratory1} for the function $f_i$ defined by on the set $Q_i \Omega_i$, with $L=n+1$. This lemma outputs a direction $\theta$ and a measure $\mu$ which we denote by $\theta_i$ and $\mu_i$ respectively. We define 
$$\Omega_{i+1} = Q_i^{-1} S_{Q_i \Omega_i, \theta _i} .$$
Equation \eqref{eq:mu0} yields that for a universal constant $c>0$,
\begin{equation} \label{eq:mu0-2}
\mu_i \left( \left\{x - x_0 : |f(x) - g(x)| > \frac{c}{n^{7.5} \log(1+n/ \eps) } \max(\epsilon, f(x)) \right\} \right) > \frac{c}{n}
\end{equation}
for all functions $g(x)$ such that $g(\alpha) < -\eps$, whenever $\alpha \in \Omega_i \setminus \Omega_{i+1}$.

Fix a constant $c'>0$ whose value will be assigned later on. Define $\delta = \tfrac{c' \eps}{16 n^{10}}$ and let
$$
N = \min \{i: ~ \exists \theta \in \Sph \mbox{ such that } |\langle x, \theta \rangle| < \delta, ~ \forall x \in \Omega_i\}.
$$
In other words, $N$ is the smallest value of $i$ such that $\Omega_i$ is contained in a slab of width $2 \delta$. Our next goal is to give an upper bound for the value of $N$.  To this end, we claim that
\begin{equation} \label{eq:voldecrease}
\Vol (\Omega_{i+1}) \leq \frac{1}{2} \Vol (\Omega_i) ,
\end{equation}
which equivalently says
$$\Vol(S_{Q_i \Omega_i, \theta_i}) \leq \frac{1}{2} \Vol(Q_i \Omega_i) .$$
Let $X \sim \mu_{Q_i \Omega_i}$ and observe that $\P(|\langle X, \theta_i \rangle| \leq 1/4) = \Vol(S_{\Omega_i, \theta_i})  / \Vol(Q_i \Omega_i)$. Clearly $\langle X, \theta_i \rangle$ is a log-concave random variable, and using that $\Cov(\Proj_{L_i} \mu_{Q_i \Omega_i}) = \Proj_{L_i}$ together with the fact that $\theta_i \in L_i \cap \mS^{n-1}$ one also has that $\langle X, \theta_i \rangle$ has variance $1$. Using that the density of a log-concave distribution of unit variance is bounded by $1$ one gets $\P(|\langle X, \theta_i \rangle| \leq 1/4) \leq 1/2$, which proves \eqref{eq:voldecrease}. It is now a simple application of Lemma \ref{lem:tech2} to see that for all $i$ there exists a direction $v_i \in \Sph$ such that 
$$
\langle v_i, \Cov(\Omega_i) v_i \rangle \leq c_1 \sqrt{n} 2^{-2i/n}.
$$
where $c_1>0$ is a universal constant. Together with Lemma \ref{lem:tech1}, this yields
$$
\diam(\Proj_{v_i} \Omega_i ) \leq 2 \sqrt{c_1} n^{5/4} 2^{-i/n}.
$$
By definition of $N$, this gives 
$$
N \leq n \log_{1/2} \tfrac{n^{5/4}}{\sqrt{c_1}} + n \log_{1/2} \delta \leq n (12 + 2c_1 + 40 \log(1+n/\eps) - \log c').
$$
Take $c'=\tfrac{\min(c,1)^2}{2^{8}(1+c_1)}$. A straightforward calculation gives
\begin{equation} \label{eq:annoying}
\frac{c}{N} > \frac{c'}{n \log(1+n/\eps)}.
\end{equation}
Finally, we define
$$
\mu(W) = \frac{1}{N} \sum_{i=1}^{N} \mu_i (W - x_0)
$$
for all measurable $W \subset \RR^n$.

For $\alpha \in \Omega \setminus \{x: |\langle x - x_0, v_N \rangle| \leq \delta \}$ consider a convex function $g(x)$ satisfying $g(\alpha) < -\eps$. Define $\tilde \alpha = \alpha - x_0$ and $\tilde g(x) = g(x + x_0) - f(x_0)$ and remark that $\tilde g(\tilde \alpha) < - \eps$. By definition of $N$, there exists $1 \leq i \leq N$ such that $\tilde \alpha \in \Omega_i \setminus \Omega_{i+1}$. Thus, equation \eqref{eq:mu0-2} gives
$$
\mu \left( \left\{x \in \Omega : |f(x) - g(x)| > \frac{c'}{2 n^{7.5} \log(1+n/ \eps) } \max(\epsilon, f(x)) \right\} \right) > \frac{c}{nN} \stackrel{\eqref{eq:annoying}}{>} \frac{c'}{n^2 \log(1+n/\eps)}.
$$
The proof is complete.

\subsection{From Lemma \ref{lem:exploratory2} to Lemma \ref{lem:exploratory1}: covering the space via regions with stable gradients} \label{sec:exploratory2}
We say that a $(z,\theta,t)$ is a {\em jolly-good triplet} if $|z| \leq \tfrac{1}{16}$ and \eqref{eq:jollygood} is satisfied for some appropriate $\delta$, namely $\delta = \frac{1}{C n^6 |\log (1+Ln/\epsilon)|}$ with $C>0$ a universal constant whose value will be decided upon later on. Intuitively given Lemma \ref{lem:exploratory2} it is enough to find a polynomial (in $n$) number of jolly-good triplets for which the corresponding set of $\theta$-directions partially covers the sphere $\mS^{n-1}$.  The notion of covering we use is the following: For a subset $H \subset \Sph$ and for $\gamma > 0$, we say that $H$ is a \emph{$\gamma$-cover} if for all $x \in \mS^{n-1}$, there exists $\theta \in H$ such that $\langle \theta, x \rangle \geq - \gamma$.

Next we explain how to find jolly-good triplets in Section \ref{sec:contraction}, and then how to find a $\gamma$-cover with such triplets in Section \ref{sec:concludingwcontraction}.

\subsubsection{A contraction lemma} \label{sec:contraction}
The following result shows that jolly-good triplets always exist, or in other words that a convex function always has a relatively big set on which the gradient map is approximately constant. Quite naturally the proof is based on a smoothing argument together with a Poincar\'e inequality.

\begin{lemma} \label{lem:contract}
Let $r, \eta, L> 0$ and $0 < \xi < 1$ such that $L>2 \eta r$. Let $\Omega \subset \R^n$ be a convex set, and $f : \Omega \rightarrow \R$ be $L$-Lipschitz and $\eta$-strongly convex, that is
$$
\nabla^2 f(x) \succeq \eta \Id,~~ \forall x \in \Omega.
$$
Let $x_0 \in \Omega$ such that $\mB(x_0, r) \subset \Omega$. Then there exist a triplet $(z, \theta, t) \in \mB(x_0, r) \times \mS^{n-1} \times [\eta r /2, +\infty)$ such that
\begin{equation}  \label{eq:funtriplet}
\mu_{\mB(z, \delta)} \left( (\nabla f)^{-1} \left( \mB\left(t \theta, \xi t\right)  \right) \right) \geq 1/2
\end{equation}
for $\delta = \frac{\xi r}{16 n^2 \log \tfrac{L}{\eta r}}$.
\end{lemma}

\begin{proof}
We consider the convolution $g = f \star h$, where $h$ is defined by
$$
h(x) = \frac{\mathbf{1}_{\{x \in \mB(0, \delta) \}}}{\Vol(\mB(0,\delta))} .
$$ 
We clearly have that $g$ is also $\eta$-strongly convex. Let $x_{min}$ be the point where $g$ attains its minimum in $\Omega$. We claim that
\begin{equation} \label{eq:gradlarge}
|\nabla g(x)| \geq \eta r / 2, ~~ \forall x \in \Omega \setminus \mB(x_{min}, r/2).
\end{equation}
Indeed by strong-convexity of $g$ we have for all $y \in \Omega$, 
$$
|\nabla g(y)| \geq \frac{1}{|y-x_{min}|} \langle \nabla g(y), y-x_{min} \rangle \geq |y-x_{min}| \eta.
$$
which proves \eqref{eq:gradlarge}. \\

Next, define $B_0 = \mB(x_0, r)$ and $D = B_0 \setminus \mB(x_{min}, r/2)$. It is clear that $\frac{\Vol(D)}{\Vol(B_0)} \geq \tfrac{1}{2}$. Let $\nu$ be the push forward of $\mu_{D}$ under $x \mapsto |\nabla g(x)|$. According to \eqref{eq:gradlarge} and by the assumption that $f$ is $L$-Lipschitz, we know that $\nu$ is supported on $[\eta r / 2, L]$. Thus, there exists some $t\in[\eta r / 2,L]$ such that $\nu([t,2t]) \geq \left (2 \log \tfrac{L}{\eta r} \right )^{-1}$. Define
$$
A = \{x \in B_0 : | \nabla g(x)| \in [t,2t] \},
$$
so we know that 
$$
\frac{\Vol(A)}{\Vol(B_0)} \geq \frac{\Vol(A)}{\Vol(D)} \frac{\Vol(D)}{\Vol(B_0)} \geq \frac{1}{4 \log  \tfrac{L}{\eta r} }.
$$
Recall that $\frac{\Vol_{n-1} (\partial B(0,r))}{\Vol_n(B(0,r))} = \frac{n+1}{r}$. Using Lemma \ref{lem:Gauss}, we now have that
$$
\frac{1}{\Vol(A)} \int_A \Delta g(x) dx \leq t \frac{\Vol_{n-1}(\partial B_0)}{\Vol(A)} = t \frac{\Vol_{n-1}(\partial B_0)}{\Vol(B_0)} \frac{\Vol(B_0)}{\Vol(A)} \leq 8 n t r^{-1} \log \tfrac{L}{\eta r}.
$$
Consequently, there exists a point $z \in A$ for which $|\nabla g(z)| \geq t$ and $\Delta g(z) \leq 8 n t r^{-1} \log \tfrac{L}{\eta r}$. In other words, by the definition of $g$, we have that
$$
\frac{1}{\Vol(\mB(z,\delta))} \int_{\mB(z, \delta) } \Delta f(x) dx \leq 8 n t r^{-1} \log \tfrac{L}{\eta r}.
$$
Fix $1 \leq i \leq n$, and define $w(x) = \langle \nabla f(x) - \nabla g(z), e_i \rangle$, where $e_i$ is the $i$-th vector of the standard basis. Note that 
$$
|\nabla w(x)| = \vert \nabla^2 f(x) e_i \vert \leq \Delta f(x).
$$
Recall that the Poincar\'e inequality for a ball (see e.g., \cite{AD03}) implies that
$$\int_{\mB(z, \delta) } |w(x)| dx \leq \delta \int_{\mB(z, \delta) } |\nabla w(x)| dx .$$
Thus combining the last three displays, and using that $\delta = \frac{\xi r}{16 n^2 \log \tfrac{L}{\eta r}}$, one obtains
$$
\frac{1}{\Vol(\mB(z,\delta))} \int_{\mB(z, \delta) } |w(x)| dx \leq  8 \delta n t r^{-1} \log \tfrac{L}{\eta r} \leq \tfrac{\xi t}{2 n}.
$$
By using the fact that $| \nabla f(x) - \nabla g(z)| \leq \sum_{i=1}^n |\langle \nabla f(x) - \nabla g(z), e_i \rangle|$ , this yields
$$
\frac{1}{\Vol(\mB(z,\delta))} \int_{\mB(z, \delta) } |\nabla f(x) - \nabla g(z)| dx \leq \xi t / 4 \leq \xi |\nabla g(z)| / 2.
$$
Finally applying Markov's inequality one obtains \eqref{eq:funtriplet} for the triplet $(z, \frac{\nabla g(z)}{|\nabla g(z)|}, |\nabla g(z)|)$.
\end{proof}

\subsubsection{Concluding the proof with the contraction lemma} \label{sec:concludingwcontraction}
We first fix some $\eta > 0$ and, at this point, suppose that $\nabla^2 f(x) \succeq \eta$ for all $x \in \Omega$. Later on we will argue that this assumption can be removed. Define $h_\Omega(x) = \sup_{y \in \Omega} \langle x, y \rangle$, the support function of $\Omega$. Consider the set
$$
\Theta = \left \{ \theta \in \Sph: ~h_\Omega(\theta) \leq \tfrac 1 8 \right \}
$$
and let $H$ be set of directions obtained from jolly-good triplets, more precisely,
$$H = \left\{\theta \in \mS^{n-1} : \exists z \in \RR^n, t \in (0, 1) \;\; \text{such that} \;\; \eqref{eq:jollygood} \;\; \text{is true with} \;\; \delta= \frac{1}{2^{28} n^6 \log (1+L n/\eta)}\right\} .$$

Define $\gamma = \tfrac{1}{16 n}$. Next, we show that $H \cup \Theta$ is a $\gamma$-cover. Let $\phi \in \mS^{n-1}$. Our objective is to find $\theta \in H \cup \Theta$ such that $\langle \theta, \phi \rangle \geq -\gamma$.

First suppose that $\phi \notin 8 \Omega$. In that case, by Hahn-Banach and since $0 \in \Omega$, there exists $w \in \RR^n$ such that $\langle \phi, w \rangle = 1$ and $\langle w, y \rangle \leq \tfrac 1 8$ for all $y \in \Omega$. In other words, we have for $\theta = \tfrac{w}{|w|}$ that
$$
h_\Omega \left (\theta \right ) \leq \frac{1}{8 |w|} \leq \frac 1 8 ,
$$
which implies that $\theta \in \Theta$. Since $\left \langle \phi, \tfrac{w}{|w|} \right \rangle \geq 0$, we are done.

We may therefore assume that $\phi/8 \in \Omega$. Since $\Cov(\mu_\Omega)=\Id$, then by Lemma \ref{lem:tech1} there exists a point $w \in \RR^n$ such that $|w| \leq n+1$ and $\mB(w,1) \subset \Omega$. Define $r = \tfrac{1}{2^{13} n^2}$ and take
$$
B_0 =  \mB(\phi / 32 + r w, r ).
$$
Note that by convexity and by the fact that $0 \in \Omega$, we have that $B_0 \subset \Omega$. We now use Lemma \ref{lem:contract} for the ball $B_0$ with $\xi=\tfrac{1}{2^{11} n^2}$, and $\delta = \frac{1}{2^{28} n^6 \log (1+L n/\eta)}$ to obtain a jolly-good triplet $(z(\theta), \theta, t)$. Denote $z=z(\theta)$. We want to show that $\langle \theta, \phi \rangle \geq - \gamma$. Observe that by convexity of $f$ and since $f$ attains its minimum at $x=0$, one has $\langle \nabla f(x), x \rangle \geq 0$ for any $x$. Thus, by definition of a jolly-good triplet one can easily see that $\langle \theta, z \rangle \geq - (\xi + \delta)$. Also by definition $z$ is in $B_0$ and thus $|32 z - \phi - 32 r w| \leq 32 r$. This implies:

\begin{align*}
\langle \theta, \phi \rangle ~& = \langle \theta, \phi - 32 z + 32 r w \rangle + 32 \langle \theta, z \rangle - 32 r \langle \theta, w \rangle \\
~& \geq - |\phi - 32 z + 32 r w| - 32 r |w| - 32 \xi - 32 \delta \geq - \tfrac{1}{16n}.
\end{align*}
This concludes the proof that $H \cup \Theta$ is a $\gamma$-cover.

Next we use Lemma \ref{lem:Caratheodory} to extract a subset $H' \subset H$ such that $|H'| \leq n+1$ and $H' \cup \Theta$ is also a $\gamma$-cover for $\mS^{n-1}$.  An application of Lemma \ref{lem:slab} with $M=2n$ now gives that there exists $v \in \Sph$ such that
$$
\Omega \cap \left (\bigcap_{\theta \in H' \cup \Theta} \left \{ x: ~ \langle x, \theta \rangle \leq \tfrac{1}{8} \right \} \right ) = \Omega \cap \left (\bigcap_{\theta \in H'} \left \{ x: ~ \langle x, \theta \rangle \leq \tfrac{1}{8} \right \} \right ) \subset S_{\Omega, v}.
$$
Finally, an application of Lemma \ref{lem:exploratory2} gives us that for all $\alpha \in \Omega \setminus S_{\Omega, v}$ and every function $g$ such that $g(\alpha) < -\eps$ one has for some $\theta \in H'$,
$$
\mu_{\mB(z(\theta), \delta)} \left( \left\{x \in \Omega : |f(x) - g(x)| > \frac{\delta}{2^{13} M \sqrt{n}} \max(\epsilon, f(x)) \right\} \right) >  \frac{1}{8} .
$$
Defining  $\mu = \frac{1}{|H'|} \sum_{\theta \in H'} \mu_{\mB(z(\theta), \delta)}$, we get
\begin{equation} \label{eq:finaldisplem2}
\mu \left( \left\{x \in \Omega : |f(x) - g(x)| > \frac{1}{2^{42} n^{7.5} \log (1+Ln/\eta)} \max(\epsilon, f(x)) \right\} \right) >  \frac{1}{16 n} .
\end{equation}
It remains to remove the uniform convexity assumption. This is done by considering the function
$$
x \mapsto f(x) + \eta |x|^2
$$
in place of $f$ in the above argument. Since $|x| \leq M \leq 2n$ for all $x \in \Omega$, the equation \eqref{eq:finaldisplem2} becomes
$$
\mu \left( \left\{x \in \Omega : |f(x) - g(x)| > \frac{c}{2^{42} n^{7.5} \log (1+Ln/\eta)} \max(\epsilon, f(x)) - 4n^2 \eta \right\} \right) >  \frac{1}{16 n}.
$$
Finally choosing $\eta = \left( \frac{\epsilon}{2^{20} n^{10}} \right)^2$ one easily obtains
$$
\mu \left( \left\{x \in \Omega : |f(x) - g(x)| > \frac{1}{2^{50} n^{7.5} \log (1+ n/ \epsilon)} \max(\epsilon, f(x)) \right\} \right) >  \frac{1}{16 n} ,
$$
which concludes the proof.

\subsection{Proof of Lemma \ref{lem:exploratory2}} \label{sec:exploratory3}

The main ingredient of the proof is the following technical result.

\begin{lemma} \label{lem:needles}
Let $\Omega \subset \RR^n$ be a domain satisfying $\mathrm{Diam}(\Omega) \leq M$. Let $f: \Omega \to [0, \infty)$ be a non-negative convex function let $g: \Omega \to \RR$ be a convex function satisfying $g(\alpha) < -\eps$, for some $\alpha \in \Omega$. Let $z \in \RR^n$ and consider the ball $B=B(z,\delta)$. Let $D \subset B$ be a set satisfying
\begin{equation} \label{eq:gradientDirection}
\left \langle \nabla f(x), \alpha - x \right \rangle \geq 0, ~~ \forall x \in D.
\end{equation}
Assume also that $\mu_B(D) \geq \frac{1}{2}$ and that $|z - \alpha| \geq n \delta$. Define 
$$A = \Bigl \{x : |f(x) - g(x)| > \frac{\delta}{2^{13} M \sqrt{n}} \max(\eps, f(x)) \Bigr \}.$$ 
Then one has $\mu_D(A) \geq 1/4$.
\end{lemma}

\begin{proof}
For $x \in \Omega$, define $\Theta_\alpha(x) = \frac{x-\alpha}{|x-\alpha|}$ and for $\theta \in \Sph$ write $N(\theta) = \Theta_\alpha^{-1} (\theta)$. Denote by $\lambda_\theta$ the one-dimensional Lebesgue measure on the needle $N(\theta)$. 
Let $\sigma_B, \sigma_D$ be the push-forward of $\mu_B, \mu_D$ under $\Theta_\alpha$. Moreover, for every $\theta \in \Sph$, the disintegration theorem ensures the existence of a probability measure $\mu_{D,\theta}$ on $N(\theta)$, defined so that for every measurable test function $h$ one has
\begin{equation} \label{eq:polar1}
\int h(x) d \mu_D(x) = \int_{\Sph} \int_{N(\theta)} h(x) d \mu_{D,\theta}(x) d \sigma_D (\theta)  
\end{equation}
(in other words, $\mu_{D,\theta}$ is the normalized restriction of $\mu_D$ to $N(\theta)$). Define the measures $\bigl (\mu_{B,\theta} \bigr)_\theta$ in the same manner.

It is easy to verify that $\sigma_D$ is absolutely continuous with respect the the uniform measure on $\Sph$, which we denote by $\sigma$. Denote $q(\theta) := \tfrac{d \sigma_D}{d\sigma} (\theta)$ and $w(\theta) := \tfrac{d \sigma_B}{d\sigma} (\theta)$.

Using Lemma \ref{lem:needles} we obtain that
\begin{equation} \label{eqLebesgue}
\frac{d \mu_{D,\theta}}{ d \lambda_\theta} (x) = \frac{\zeta_n}{ \Vol(D) q(\theta) } |x-\alpha|^{n-1} \mathbf{1}_{ \{x \in D\}} ,
\end{equation}
and
\begin{equation} \label{eq:Lebesgue2}
\frac{d \mu_{B,\theta}}{ d \lambda_\theta} (x) = \frac{\zeta_n} {\Vol(B) w(\theta) } |x-\alpha|^{n-1} \mathbf{1}_{ \{x \in B\}} ,
\end{equation}
where $\zeta_n$ is a constant depending only on $n$.

For every $\theta \in \Sph$, define $L(\theta)$ to be the length of the interval $N(\theta) \cap B$. Consider the set 
$$
\mathcal{L} = \left \{ \theta : ~ L(\theta) > \frac{\delta}{32 \sqrt{n}} \right \} .
$$
According to Lemma \ref{lem:ball-needle} we have that 
$$
\int_{\Sph \setminus \mathcal{L}} w(\theta) d \sigma(\theta) \leq \frac{1}{8}.
$$
Now, since $D \subset B$ and $\mu_B(D) \geq \tfrac 1 2$, we have that $q(\theta) \leq 2 w(\theta)$ for all $\theta \in \Sph$, which gives
$$
\sigma_D(\mathcal{L}) = \int_{\mathcal{L}} q(\theta) d \sigma(\theta) \geq \frac{3}{4}.
$$
Next, consider the set
$$
\mathcal{S} = \left \{ \theta \in \Sph; ~ q(\theta) \geq \frac{w(\theta)}{4} \right \}.
$$
Since $\int_{\Sph} \frac{q(\theta)}{w(\theta)} d \sigma_B(\theta) = 1$ we have
$$
\sigma_D(\mathcal{S}) = \int_{\mathcal{S}} \frac{q(\theta)}{w(\theta)} d \sigma_B(\theta) = 1 - \int_{\Sph \setminus \mathcal{S}} \frac{q(\theta)}{w(\theta)} d \sigma_B(\theta) \geq \frac{3}{4}.
$$
Using a union bound, we have that $\sigma_D(\mathcal{L} \cap \mathcal{S}) \geq \frac{1}{2}$.

Fix $\theta \in \mathcal{L} \cap \mathcal{S}$, we would like to give a lower bound on $\mu_{D,\theta} (A)$. In view of Lemma \ref{lem:onedim}, we thus need an upper bound on the density of $\mu_{D,\theta}$. Recall that $\theta \in \mathcal{S}$, implies $\frac{q(\theta)}{w(\theta)} \geq \tfrac{1}{4}$ and that by \eqref{eqLebesgue} and \eqref{eq:Lebesgue2}, we have for all $x \in N(\theta) \cap B$,
\begin{equation} \label{eq:compqw}
\frac{d \mu_{D,\theta}}{ d \mu_{B,\theta}} (x) = \frac{\Vol(B) w(\theta)}{\Vol(D) q(\theta)} \mathbf{1}_{x \in D} \leq 8. 
\end{equation}
Denote $[a,b] = B \cap N(\theta)$ for $a,b \in \RR^n$. Assume that $a$ is the interior of the interval $[\alpha,b]$ (if this is not the case, we simply interchange between $a$ and $b$). By the assumption $\theta \in \mathcal{L}$, we know that $|b-a| \geq \frac{\delta}{32 \sqrt{n}}$. Writing $Z = \frac{\zeta_n }{\Vol(B) w(\theta) }$ so that, according to \eqref{eq:Lebesgue2},
$$
\frac{d \mu_{B,\theta}}{ d \lambda_\theta} (x) = Z |x-\alpha|^{n-1} \mathbf{1}_{ \{x \in B\}},
$$
and since $\mu_{B,\theta}$ is a probability measure,
$$
Z^{-1} = \int_a^b |x-\alpha|^{n-1} dx
$$
where, by slight abuse of notation we assume that $a,b,\alpha \in \RR$. Thus,
$$
Z \leq \frac{32 \sqrt{n}}{\delta |a-\alpha|^{n-1}}.
$$
Combined with \eqref{eq:compqw}, this finally gives
\begin{align*} 
\frac{d \mu_{D,\theta}}{ d \lambda_\theta} (x) ~& \leq 2^8 \frac{\sqrt{n}}{\delta} \frac{|x-\alpha|^{n-1}}{ |a-\alpha|^{n-1}} \leq 2^8 \frac{\sqrt{n}}{\delta}  \left (\frac{|b-\alpha|}{|a-\alpha|}\right )^{n-1} \\
~& = 2^8 \frac{\sqrt{n}}{\delta}  \left (1 + \frac{|b-a|}{|a - \alpha|} \right )^{n-1} \leq 2^8 \frac{\sqrt{n}}{\delta}  \left (1 + \frac{2\delta}{n \delta - \delta} \right )^{n-1} \leq 2^8 e^2 \frac{\sqrt{n}}{\delta} ,
\end{align*}
where in the second to last inequality we used the assumption that $|z-\alpha| \geq n \delta$. 

Define the map $U: \RR \to N(\theta)$ by
$$
U(x) = \alpha + M (|\alpha| - x ) \theta
$$
and consider the functions $\tilde f(x) = f(U(x))$ and $\tilde g(x) = g(U(x))$. Denote $x_0 = \min U^{-1}( D \cap N(\theta))$ and remark that $x_0\in [|\alpha| -1, |\alpha| ]$. Note that, thanks to equation \eqref{eq:gradientDirection}, the assumption \eqref{eq:assumpDerive} holds for the functions $\tilde f, \tilde g$ and the points $x_0, |\alpha|$. We can now invoke Lemma \ref{lem:onedim} for these functions with $\mu$ being the pullback of $\mu_{D,\theta}$ by $U(x)$. According to the above inequality one may take $\beta = 2^8 e^2 \frac{M\sqrt{n}}{\delta}$ and obtain
$$
\mu_{D,\theta} (A) \geq \frac{1}{2}.
$$
Integrating over $\theta \in \mathcal{L} \cap \mathcal{S}$ concludes the proof:
$$
\mu_D(A) \geq \int_{\mathcal{S} \cap \mathcal{L}} \mu_{D,\theta} (A) d \sigma_D(\theta) \geq \frac{1}{2} \sigma_D(\mathcal{L} \cap \mathcal{S} ) \geq \frac{1}{4}.
$$
\end{proof}

\begin{proof} [Proof of Lemma \ref{lem:exploratory2}]
Suppose that $(z, \theta, t)$ satisfy equation \eqref{eq:funtriplet}. Fix $\alpha \in \Omega$ satisfying $\langle \alpha, \theta \rangle \geq \tfrac{1}{8}$ and a function $g(x)$ satisfying $g(\alpha) < -\eps$. 
Define $B = B(z, \delta)$ and $D = \left \{x \in B; |\nabla f(x) - \theta t| < \frac{1}{16} n^{-2} t \right \}$. Let $\mu_B$ be the uniform measure on $B$. According to \eqref{eq:funtriplet}, we know that $\mu_B(D) \geq \tfrac{1}{2}$. Now, for all $x \in D$ we have that $\nabla f(x) = t(\theta + y)$ with $|y| < \tfrac{1}{16} n^{-2}$ so we get
\begin{align*} 
\left \langle \nabla f(x), \frac{\alpha - x}{|\alpha - x|} \right \rangle &~= \frac {t} {|\alpha - x|} \left (\langle \alpha, \theta \rangle + \langle \alpha - x, y \rangle - \langle x, \theta \rangle \right ) \\
&~> \frac{t}{|\alpha - x|} \left (\tfrac 1 8 - \tfrac{1}{16}(|\alpha| + |x|)n^{-2}- |x| \right ) \geq 0, ~~ \forall x \in D
\end{align*}
where we used the fact that $D \subset B$ and so $|x|< |z| + \delta \leq \tfrac 1 {16}$ and the fact that $|\alpha| \leq 2n$. Note that the above implies the assumption \eqref{eq:gradientDirection}. Moreover remark that 
$$
|z-\alpha| \geq \tfrac 1 4 - \tfrac 1 8 \geq \tfrac 1 8 \geq n \delta.
$$ 
We can thus now invoke Lemma \ref{lem:needles} to get $\mu_B(A) \geq 1/8$ where
$$
A = \bigl \{x \in \Omega : ~ |f(x) - g(x)| > \tfrac{\delta}{2^{13} M \sqrt{n}} \max(\eps, f(x)) \bigr \}.
$$
This completes proof. 
\end{proof}

\subsection{Technical lemmas} \label{sec:exploratory4}
We gather here various technical lemmas.

\begin{lemma} \label{lem:tech1}
Let $C$ be a convex body in $\RR^n$. Then
\begin{equation} \label{eq:diamisotropic}
\diam(C) \leq (n+1) \| \Cov(\mu_C) \|^{1/2}. 
\end{equation}
On the other hand, if $\Cov(\mu_C) \succeq \Id$ then $C$ contains a ball of radius $1$.

Furthermore, for all $v \in \Sph$ one has
$$
\sup_{x \in C} \langle v, x \rangle - \inf_{x \in C} \langle v, x \rangle \leq (n+1) \langle v, \Cov(\mu_C), v \rangle^{1/2}.
$$
\end{lemma}

\begin{proof}
The first and second parts of the Lemma are found in \cite[Section 3.2.1]{MR3185453}. For the second part, we write $C' = \Cov(C)^{-1/2} C$ and $u = \tfrac{\Cov(C)^{1/2} v}{|\Cov(C)^{1/2} v|}$. We have
\begin{align*}
\sup_{x \in C} \langle v, x \rangle - \inf_{x \in C} \langle v, x \rangle ~& = \sup_{x \in C'} \langle v, \Cov(C)^{1/2} x \rangle - \inf_{x \in C'} \langle v, \Cov(C)^{1/2} x \rangle  \\
&= \sup_{x \in C'} \langle \Cov(C)^{1/2} v, x \rangle - \inf_{x \in C'} \langle \Cov(C)^{1/2} v,  x \rangle\\
& = |\Cov(C)^{1/2} v| \left (\sup_{x \in C'} \langle u, x \rangle - \inf_{x \in C'} \langle u, x \rangle \right ) \stackrel{\eqref{eq:diamisotropic}}{\leq} (n+1) |\Cov(C)^{1/2} v|.
\end{align*}
\end{proof}

\begin{lemma} \label{lem:tech2}
Let $C \subset D \subset \RR^n$ be two convex bodies with $0 \in C$. Suppose that $\frac{\Vol(C)}{\Vol(D)} \leq \delta$, then there exists $u \in \mS^{n-1}$ such that
\begin{equation} \label{eq:covdecr}
\langle u, \mathrm{Cov}(\mu_{C}) u \rangle \leq c \sqrt{n} \delta^{2/n} \langle u, \mathrm{Cov}(\mu_{D}) u \rangle.
\end{equation}
where $c>0$ is a universal constant.
\end{lemma}
\begin{proof}
Define $\mu = \mu_{D}$ and $\nu=\mu_{C}$. By applying a linear transformation to both $\mu$ and $\nu$, we can clearly assume that $\Cov(\mu) = \mathrm{Id}$. Let $f(x)$ be a log-concave probability density in $\RR^n$. According to \cite[Corollary 1.2 and Lemma 2.7]{MR2276540}, we have that
\begin{equation} \label{eq:Klartag}
c_1 \leq \left (\sup_{x \in \RR^n} f(x) \right )^{1/n} \left (\det \Cov (f)\right )^{1/2n} \leq c_2 n^{1/4}
\end{equation}
where $c_1,c_2>0$ are universal constants. Denote by $f(x)$ and $g(x)$ the densities of $\mu$ and $\nu$, respectively. Since $\mu,\nu$ are indicators, we have that 
$$
\sup_{x \in \RR^n} f(x) = f(0) = \delta g(0) = \delta \sup_{x \in \RR^n} g(x).
$$ 
We finally get
\begin{align*}
\left (\det \Cov(\nu)\right )^{1/n} ~& \stackrel{\eqref{eq:Klartag}}{\leq} c_2^2 \sqrt n g(0)^{-2/n}  \\
& = c_2^2 \delta^{2/n} \sqrt{n} f(0)^{-2/n} \\
& \stackrel{\eqref{eq:Klartag}}{\leq} (c_2/c_1)^2 \sqrt{n} \left (\det \Cov(\mu)\right )^{1/n} \delta^{2/n} = (c_2/c_1)^2 \sqrt{n} \delta^{2/n}.
\end{align*}
The lemma follows by taking $u$ to be the eigenvector corresponding to the smallest eigenvalue of $\Cov(\nu)$.
\end{proof}

\begin{lemma} \label{lem:Gauss}
Let $g$ be a convex function defined on a Euclidean ball $B \subset \R^n$. Let $A \subset B$ be a closed set such that $\forall x \in A$, $|\nabla g(x)| \leq t$. Then
$$
\int_A \Delta g(x) dx \leq t \ \Vol_{n-1}(\partial B).
$$
\end{lemma}

\begin{proof}
Since $g$ is convex, we can write
$$
g(x) = \sup_{y \in B} w_y (x)
$$
where $w_y(x) = \langle x-y, \nabla g(y) \rangle + g(y)$. Define
$$
\tilde g(x) = \sup_{y \in A} w_y(x).
$$
Clearly $\tilde g$ is convex and $\tilde g (x) = g(x)$ for all $x \in A$. Moreover $|\nabla \tilde g(x)| \leq t$ for all $x \in \R^n$. Using Gauss's theorem, we have
$$
\int_A \Delta g(x) dx \leq \int_B \Delta \tilde g(x) dx = \int_{\partial B} \langle \nabla \tilde g(x), n(x) \rangle d \mathcal{H}_{n-1}(x) \leq t \ \Vol_{n-1}(\partial B) ,
$$
which concludes the proof
\end{proof}

Let $\gamma > 0$. Recall that we say that $H \subset \Sph$ is a \emph{$\gamma$-cover} if for all $x \in \mS^{n-1}$, there exists $\theta \in H$  satisfying
\begin{equation} \label{eq:scalarp}
\langle \theta, x \rangle \geq - \gamma.
\end{equation}

\begin{lemma} \label{lem:Caratheodory}
Let $H \subset \Sph$ be a $\gamma$-cover. Then there exists a subset $I \subset H$ with $|I| \leq n+1$ such that $I$ is a $\gamma$-cover.
\end{lemma}
\begin{proof}
We first claim that there is a point $y \in \mathrm{Conv} (H)$ with $|y| \leq \gamma$. Indeed, if we assume otherwise then by Hahn-Banach there exists $\tilde \theta \in \Sph$ such that $\langle \theta, \tilde \theta \rangle > \gamma$ for all $\theta \in H$, which means the vector $-\tilde \theta$ violates the assumption \eqref{eq:scalarp}. By Caratheodory's theorem, there exists $I \subset H$ with $|I| \leq n+1$ such that $y \in \mathrm{Conv} (I)$. Write $I = (\theta_1,...,\theta_{n+1})$. Now let $x \in \RR^n$ with $|x| \leq 1$. Then since $\langle x, y \rangle \geq - \gamma$, we have
$$
\sum_{i=1}^{n+1} \alpha_i \left \langle x, \theta_i \right \rangle \geq -\gamma  
$$
for some non-negative coefficients $\{\alpha_i\}_{i=1}^{n+1}$ satisfying $\sum_{i = 1}^{n+1} \alpha_i = 1$. Thus there exists $\theta \in I$ for which \eqref{eq:scalarp} holds.
\end{proof}

\begin{lemma} \label{lem:slab}
Let $\Omega \subset \R^n$ be a convex set with $\diam(\Omega) \leq M$ and such that $0 \in \Omega$. Let $H$ be a $\gamma$-cover. Then there exists $\tilde{\theta} \in \mS^{n-1}$ such that
$$\{ \alpha \in \Omega : \forall \theta \in H, \langle \alpha, \theta \rangle < M \gamma \} \subset \{\alpha \in \Omega : |\langle \alpha, \tilde \theta \rangle| \leq 2 M \gamma \} .$$
\end{lemma}

\begin{proof}
Since $\{ \alpha \in \Omega : \forall \theta \in H, \langle \alpha, \theta \rangle < M \gamma \}$ is a convex set which contains $0$, showing that it does not contain a ball of radius $2 M \gamma$ is enough to show that it is included in some slab $\{\alpha \in \Omega : |\langle \alpha, \tilde \theta \rangle| \leq 2 M \gamma \}$. Now suppose that our set of interest $\{ \alpha : \forall \theta \in H, \langle \alpha, \theta \rangle < M \gamma \}$ actually contains a ball $\mB(x, 2 M \gamma)$ with $|x| \in (0, M)$. Let $\theta \in H$ be such that $\langle \frac{x}{|x|}, \theta \rangle \geq - \gamma$, and thus in particular $\langle x, \theta \rangle \geq - M \gamma$. Then one has by the inclusion assumption that $\langle \theta, x + 2 M \gamma \theta \rangle < M \gamma$, but on the other hand one also has $\langle \theta, x + 2 \gamma M \theta \rangle \geq \gamma M$ which yields a contradiction, thus concluding the proof.
\end{proof}

\begin{lemma} \label{lem:ball-needle}
Let $\delta > 0$, $x_0 \in \RR^n$, $B=\mB(x_0,\delta)$ and $\alpha \in \RR^n \setminus B$ . For $x \in \RR^n$, define $\Theta_\alpha(x) = \frac{x-\alpha}{|x-\alpha|}$, and let $\sigma_B$ be the push-forward of $\mu_B$ under $\Theta_\alpha$. For every $\theta \in \Sph$, define $L(\theta)$ to be the length of the interval $\Theta_\alpha^{-1}(\theta) \cap B$. Then one has
$$
\sigma_B \left ( \theta : ~ L(\theta) > \frac{\delta}{32 \sqrt{n}} \right ) \geq \frac{7}{8}.
$$
\end{lemma}

\begin{proof}
Note that, by definition,
$$
x \in B \mbox{ and } x + \frac{\delta}{32 \sqrt{n}} \frac{\alpha - x}{|\alpha - x|} \in B \Rightarrow  L\left(\Theta_\alpha(x)\right) > \frac{\delta}{32 \sqrt{n}} .
$$
Furthermore it is easy to show that for all $y \in B$,
$$
y + \frac{\delta}{32 \sqrt{n}} \frac{\alpha - x_0}{|\alpha - x_0|} \in B \Rightarrow y + \frac{\delta}{32 \sqrt{n}} \frac{\alpha - y}{|\alpha - y|} \in B.
$$
Thus letting $X \sim \mu_B$ we see that the lemma will be concluded by showing that
$$
\P \left (X + \frac{\delta}{32 \sqrt{n}} \frac{\alpha - x_0}{|\alpha - x_0|} \in B \right ) \geq \frac{7}{8}.
$$
Defining $\tilde B = \mB\left (x_0 - \frac{\delta}{32 \sqrt{n}} \frac{\alpha - x_0}{|\alpha - x_0|}, \delta \right )$, the statement boils down to proving that $\P(X \in \tilde B) \geq 7/8$. By applying an affine linear transformation to both $B$ and $\tilde B$, this is equivalent to
$$
\frac{\Vol \Bigl (\mB \left (- \frac{c}{2 \sqrt{n}} e_1 ,1\right ) \cap B \left  (\frac{c }{2 \sqrt{n}} e_1, 1 \right  ) \Bigr)}{\Vol(\mB(0,1))} \geq \frac{7}{8}
$$
where $e_1$ is the first vector of the standard basis. Next, by symmetry around the hyperplane $e_1^\perp$, we have
\begin{align*}
\frac{\Vol \Bigl (\mB \left (- \frac{1}{64 \sqrt{n}} e_1 ,1\right ) \cap \mB \left  (\frac{1}{64 \sqrt{n}} e_1, 1 \right  ) \Bigr)}{\Vol(\mB(0,1))} = \frac{2 \Vol \Bigl (\mB \left (- \frac{1}{64 \sqrt{n}} e_1 ,1\right ) \cap \left  \{ x; \langle x, e_1 \rangle \geq 0  \right \} \Bigr)}{\Vol(\mB(0,1))}.
\end{align*}
Thus, it is enough to show that $\P \left (|Z| > \tfrac{1}{64 \sqrt n} \right ) \geq \frac{7}{8}$ where $Z = \langle X', e_1 \rangle $ and $X' \sim \mu_{\mB(0,1)}$. Observe that $\Var[Z] \geq \frac{1}{8 n}$ and that $Z$ is log-concave (in particular the density of $Z/\Var[Z]$ is bounded by $1$). This implies that for any $t > 0$
$$
\P \left (|Z| < t \sqrt{\Var[Z]} \right ) < 2 t ,
$$
and thus the lemma follows by taking $t = \tfrac{1}{16}$.
\end{proof}

\begin{lemma} \label{lem:needlestech}
Let $A \subset \R^n$. For $x \in \RR^n$, define $\Theta_\alpha(x) = \frac{x-\alpha}{|x-\alpha|}$, and let $\sigma_A$ be the push-forward of $\mu_A$ under $\Theta_\alpha$. Assume that $\sigma_A$ is absolutely continuous with respect the the uniform measure $\sigma$ on $\Sph$ and denote $q(\theta) := \tfrac{d \sigma_A}{d\sigma} (\theta)$. Finally let $\mu_{A,\theta}$ be the normalized restriction of $\mu_A$ on $N(\theta) = \Theta_\alpha^{-1}(\theta)$, defined so that for every measurable test function $h$ one has
\begin{equation} \label{eq:polar1}
\int h(x) d \mu_D(x) = \int_{\Sph} \int_{N(\theta)} h(x) d \mu_{A,\theta}(x) d \sigma_A (\theta)  .
\end{equation}
Denoting $\zeta_n$ for the $(n-1)$-dimensional Hausdorff measure of $\Sph$ one then obtains
\begin{equation} \label{eqLebesgue3}
\frac{d \mu_{A,\theta}}{ d \lambda_\theta} (x) = \frac{\zeta_n}{\Vol(A) q(\theta) } |x-\alpha|^{n-1} \mathbf{1}_{ \{x \in B\}}.
\end{equation}
\end{lemma}

\begin{proof}
First observe that the existence of $\mu_{A,\theta}$ is ensured by the disintegration theorem. Now remark that using the integration by polar coordinates formula we have for every measurable test function $\varphi$,
$$
\int_{\RR^n} \varphi(x) d x = \zeta_n \int_{\Sph} \int_{0}^\infty r^{n-1} \varphi(\alpha + r \theta )  dr d \sigma(\theta) .
$$
Now, by definition of $q(\cdot)$, we have for every test function $\varphi$,
$$
\int_{\Sph} \int_{0}^\infty r^{n-1} \varphi(\alpha + r \theta )  dr d \sigma(\theta) = \int_{\Sph} \int_{0}^\infty r^{n-1} q(\theta)^{-1} \varphi(\alpha + r \theta )  dr d \sigma_A(\theta).
$$
Taking $\varphi(x) = h(x) \mathbf{1}_{x \in A}$, we finally get
\begin{align*}
\int h(x) d \mu_A(x) &~= \frac{1}{\Vol(A)} \int_{A} h(x) d x \\
&~ = \frac{\zeta_n}{\Vol(A)} \int_{\Sph} \int_{0}^\infty r^{n-1} q(\theta)^{-1} h(\alpha + r \theta) \mathbf{1}_{\{\alpha + r \theta \in A\}}  dr d \sigma_A(\theta).
\end{align*}
Since the above is true for every measurable function $h$, together with equation \eqref{eq:polar1} we get that for every function $h$ and every $\theta \in \Sph$, one must have
\begin{align*}
\int_{N(\theta)} h(x) d \mu_\theta(x) &~= \frac{\zeta_n}{Vol(D) q(\theta)} \int_{0}^\infty r^{n-1} h(\alpha + r \theta) \mathbf{1}_{\{\alpha + r \theta \in A\}}  dr \\
&~= \frac{\zeta_n}{\Vol(A) q(\theta)}  \int_{N(\theta)} |x-\alpha|^{n-1}  h(x) \mathbf{1}_{\{x \in A\}}  d\lambda_\theta(x)
\end{align*}
and the claimed identity \eqref{eqLebesgue3} follows.
\end{proof}

\bibliographystyle{plainnat}
\bibliography{newbib}
\end{document}